\documentclass[12pt,a4paper]{article}

\usepackage{amsmath,mathtools,pifont,amsfonts,amssymb,latexsym,multicol}
\usepackage{float,xspace,calc,theorem,ifthen}
\usepackage{fancyhdr,enumerate,epsfig}
\usepackage{hyperref}

\numberwithin{equation}{section}
\numberwithin{figure}{section}
\theoremstyle{change}
\newtheorem{theorem}{Theorem} [section]
\newtheorem {lemma}[theorem]{Lemma}

{\theorembodyfont{\normalfont}\newtheorem {remark}[theorem]{Remark}}
{\theorembodyfont{\normalfont}\newtheorem {remarks}[theorem]{Remarks}}
{\theorembodyfont{\normalfont}\newtheorem {example}[theorem]{Example}}
\newcommand{\beq}{\begin{equation}}
\newcommand{\eeq}{\end{equation}}
\newcommand{\Leq}[1]{\label{#1}\end{equation}}
\newcommand{\beqn}{\begin{eqnarray}}
\newcommand{\eeqn}{\end{eqnarray}}
\newcommand{\beqno}{\begin{eqnarray*}}
\newcommand{\eeqno}{\end{eqnarray*}}
\newcommand {\vep}{\varepsilon}
\newcommand {\LA}{\left\langle}
\newcommand {\RA}{\right\rangle}
\newcommand {\eh}{{\textstyle \frac{1}{2}}}

\newcommand {\sign}{{\rm sign}}
\renewcommand{\Im}{{\rm Im}}
\renewcommand{\Re}{{\rm Re}}
\newcommand {\tr}{{\rm tr}}
\newcommand {\bC}{{\mathbb C}}

\newcommand {\bN}{{\mathbb N}}
\newcommand {\bR}{{\mathbb R}}
\newcommand {\bS}{{\mathbb S}}

\newcommand {\bZ}{{\mathbb Z}}
\newcommand{\rstr}{{\upharpoonright}}
\newcommand{\bem}{\l(\! \begin{array}}
\newcommand{\eem}{\end{array}\!\ri)}

\newcommand{\qmbox}[1]{\quad\mbox{#1}\quad}
\usepackage[T1]{fontenc}
\begin{document}
\title {Regularisation by Hamiltonian extension}

\date{\today}

\author{Andreas Knauf\thanks{
Department of Mathematics,
Friedrich-Alexander-University Erlangen-N\"urnberg.
Cauerstr.~11, D-91058 Erlangen, 
Germany, \texttt{knauf@math.fau.de}}}

\maketitle
\begin{abstract}
We consider the Kepler potential and its relatives $q\mapsto -\|q\|^{-2(1-1/n)}$, $n\in\bN$
in arbitrary dimension $d$. We derive a unique real-analytic symplectic extension 
of phase space on which the Hamiltonian flow is complete and still real-analytic.
\end{abstract}
\tableofcontents
%
\section{Introduction}

There has been a long tradition to regularise the two-body problem of celestial mechanics,
in order to deal with (near) collision orbits. 
Here we present a method that works not only for any spatial dimension, but also for
a known discrete family of homogeneous potentials. 
Specifically, we consider for dimension $d\in \bN$
(negative) potentials of (negative) homogeneity 
$\alpha_n:=2(1-1/n)$ with $n\in\bN$, $n\ge2$: 
\[ U_n: \bR^d_q \backslash \{0\} \to \bR\qmbox{,} U_n(q) = Z \|q\|^{-\alpha_n} \, ,\]
with strength of interaction  $Z>0$. For reasons of comparison, 
we also look at the free case with $n=1$ and $U_1$ constant.

These potentials were discussed by {\sc McGehee} in \cite{McG}.
For them a regularisation of Levi-Civita type in dimension 
$d=2$, that is, $\bR^d_q\cong \bC$,
is of the form $q=Q^n$, $n=2$ being Kepler's case. 
Not knowing a name for them, we tentatively call them {\em McGehee potentials}.

Some typical trajectories for motion in these potentials at energy 0 are shown in 
Figure \ref{fig:Tschirnhaus}.
The solution curves for $U_3(q)= Z\|q\|^{-4/3}$ at total energy zero are Tschirnhaus 
cubics.

The unregularised\,\footnote{Unregularised quantities wear a hat.} 
phase space is the cotangent bundle 
\[ \widehat{P} := T^* \widehat{M} \quad \mbox{of configuration space }
 \widehat{M} := \bR^d_q \backslash \{0\} \]
with its canonical symplectic form $\widehat{\omega}$.
Then with mass $m > 0$ the Hamiltonian function
\beq \textstyle
\widehat{H}: \widehat{P} \to \bR \qmbox{,} 
\widehat{H} (q,p) = \frac{\|p\|^2} {2m} - U_n (q) 
\Leq{Hamiltonian}
is real-analytic.
Our main result is a generalisation of the physical
case $d=3$ in \cite[Theorem 5.1]{Kn},
to arbitrary dimensions $d\in\bN$, and additionally to all $\alpha_n$.\\ 
It is similar in spirit to the regularisation of the Kepler case
by {\sc Ligon} and {\sc Schaaf} in \cite{LS}
which was discussed by {\sc Heckmann} and {\sc de Laat} in \cite{HdL}. 
However, one important difference already for $\alpha_2=1$ is that our regularisation
includes all energies.
\begin{theorem}[completion by analytic extension] \label{thm:complete}\quad\\
For $n\in\bN$ and $d\ge 2$ 
there exists a unique Hamiltonian system $(P,\omega,H)$
which is a real-analytic extension
of $\big( \widehat{P},\widehat{\omega},\widehat{H} \big)$, 
having a complete vector field $X_H$. 
\end{theorem}
\begin{remark}[existence and uniqueness] \label{rem:eu}\quad\\
For dimension $d=1$, too there exists such a Hamiltonian system $(P,\omega,H)$,
but it is not uniquely defined (see Remark \ref{rem:one}.1).

For $d\ge2$ uniqueness up to isomorphisms is granted because every collision orbit 
must contain exactly one collision point, and there is only one {\em continuous} 
way to continue it after collision, comparing with nearby non-collision orbits.
So the article basically shows {\em existence} of such a $(P,\omega,H)$.
\hfill $\Diamond$
\end{remark}
Although the solutions of the Hamiltonian equations only have an essential 
dependence on $n$, we kept the constants $m$ and $Z$ for later applications.

Two basic ideas lead to our construction:
\begin{enumerate}[$\bullet$]
\item 
We consider (physical) time $t$ as a phase space function $T$ near collision and
complement it by suitable constants of motion to obtain a symplectic chart. 
These functions on phase space are simply total energy,
angular momentum and a generalised Laplace-Runge-Lenz vector.
\item 
We deal with the fact that -- unlike for free and Kepler motion -- for $n\ge3$ the Hamiltonian
flow is not super-integrable. This follows from Bertrand's Theorem and means that for 
negative energies the typical orbit is not closed. So a generalised 
Laplace-Runge-Lenz vector can only exist {\em locally} in phase space.
But it turns out that this suffices for our construction.
\end{enumerate}
{\bf Acknowledgement:} I sincerely thank Alain Albouy for many 
discussions about regularisation and the McGehee potentials.
%
\section{Analytic regularisation}
To get a feel for the dynamics associated to the potentials $U_n$, in the 
first subsection we shortly present their regularisation \`a la Levi-Civita.
Then in Section \ref{sub:proof} we turn to the proof of 
Theorem \ref{thm:complete}, which makes some use of 
Section \ref{sub:Levi:Civita}.
%
\subsection{Levi-Civita-like regularisation} \label{sub:Levi:Civita}
Because of conservation of angular momentum, for any dimension $d$ the 
solution curves lie in a plane.
So we assume that $\widehat{M} = \bC^\times$ 
so that $\widehat{P} \cong \bC^\times \times \bC$. 

Up to time parameterisation, the curves that correspond to Hamiltonian orbits
at energy $E$ (and do not touch the boundary circle $U_n \equiv -E$ of Hill's 
region for $E<0$) coincide with the unit speed 
geodesics for the Jacobi-Maupertuis metric $\widehat{g}_E$ on $\widehat{M}$ that
equals $W := \tfrac 2m(E+U_n)$ times the Euclidean metric. 

Although the hypersurfaces
$(\widehat{H}-E)^{-1}(0)$ and $\big((\widehat{H}-E)/W)\big)^{-1}(0)$ of $\widehat P$ coincide, the factor of time reparameterisation is non-trivial: 
\beq \textstyle
\frac{\mathrm{d}s}{\mathrm{d}t}=W(q(t)) \qmbox{or conversely}
\frac{\mathrm{d}t}{\mathrm{d}s}=\frac{1}{W(Q(s))} =\tfrac m2 \frac{|Q(s)|^{2(n-1)}}{Z + E\,|Q(s)|^{2(n-1)}}\, .
\Leq{time:rep}
When, as in \cite[Section 8] {McG}, we consider the $n$-fold branched covering 
\[\pi: \bC\to\bC \qmbox{,} \pi(Q)=Q^n\qquad\mbox{and }\
\hat{\pi} := \pi|_{\bC^\times} \, ,\]
then $\hat{\pi}$ is unbranched.  
The lift of $U_n$ has the form $U_n\circ \hat{\pi}(Q)=Z |Q|^{2(1-n)}$, so that the lift 
of the Jacobi metric $\widehat{g}_E$ equals, for $\Re(Q)=Q_1$ and  $\Im(Q)=Q_2$
\beq   
\hat{\pi}^*\, \widehat{g}_E(Q) = 
2n^2m \big(Z+ E (Q\bar{Q})^{2n-2}\big) (dQ_1 \otimes dQ_1 + dQ_2 \otimes dQ_2) \, . 
\Leq{lifted:metric}
Since $n\in \bN$, this uniquely extends to 
a real-analytic metric $\mathbf{g}_E$ on $\bC$ if $E\ge 0$, respectively
on lifted Hill's region $\{Q\in \bC \mid  |Q|^{2(n-1)} < - Z/E\}$ if $E < 0$.\\
So the geodesic flow on the covering surface with metric $\mathbf{g}_E$
is real-analytic in the time parameter $s$ and the initial conditions.\\ 
The factor $\frac{\mathrm{d}t}{\mathrm{d}s}$ of time reparameterisation is real-analytic too, see \eqref{time:rep}. 

\begin{figure}[h]
\centerline{
\includegraphics[height=40mm]{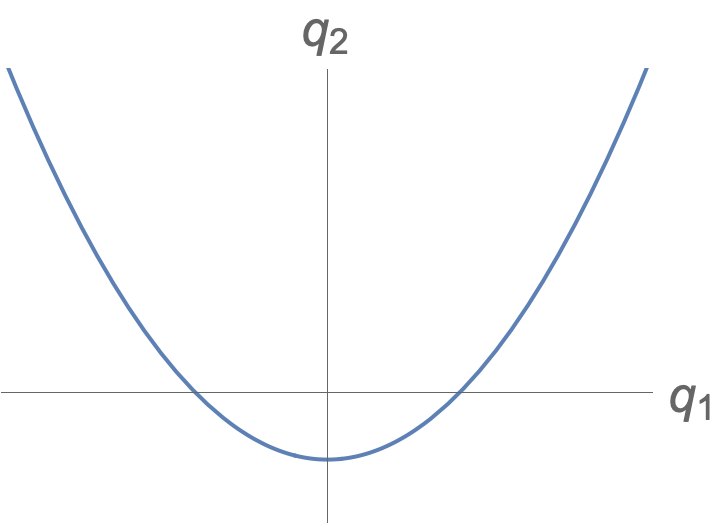}
\hfill
\includegraphics[height=40mm]{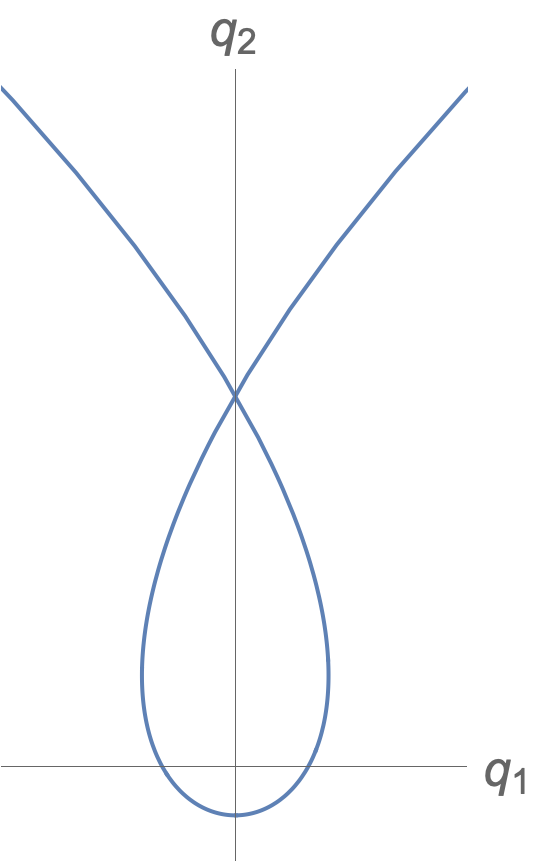}
\hfill
\includegraphics[height=40mm]{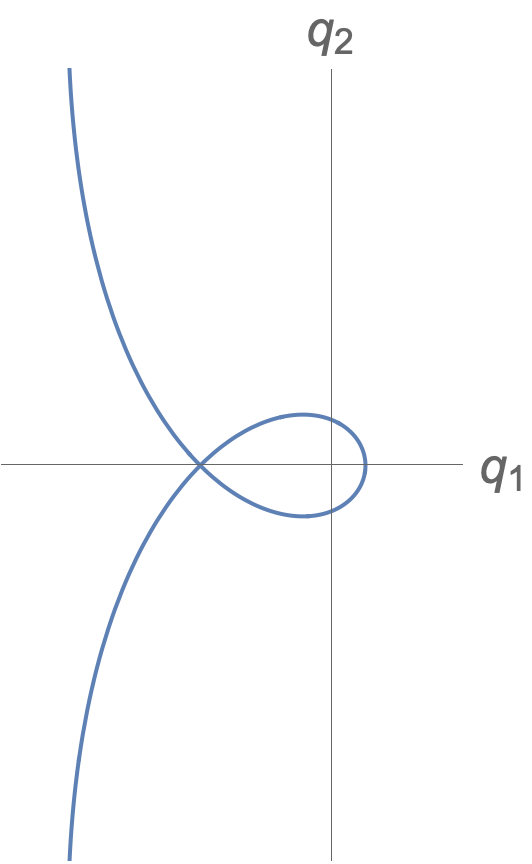}
\hfill
\includegraphics[height=40mm]{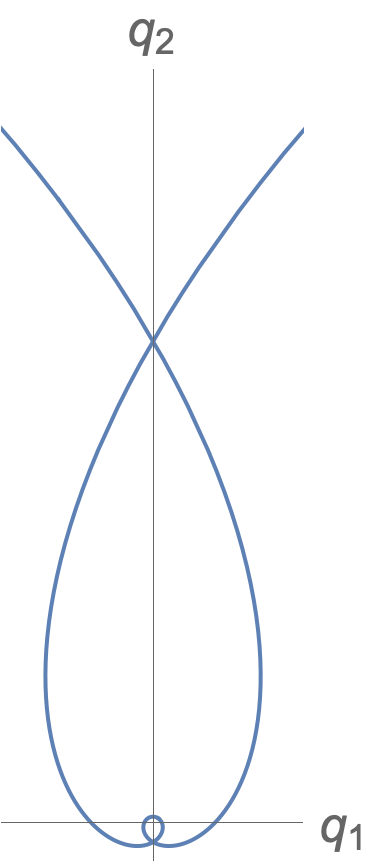}
}
\caption{Solution curves with potential $U_n$ for energy zero. 
From left to right:
$n=2$ (Kepler parabola), $n=3$ (Tschirnhaus curve), $n=4$, $n=6$. Despite 
their appearance (due to the small parameter range shown),
their asymptotic directions $\lim_{t\to\pm\infty}\frac{q(t)}{|q(t)|}$ 
coincide if $n$ is even, and are opposite if $n$ is odd.}
\label{fig:Tschirnhaus}
\end{figure}

In particular, for $E=0$, the geodesics are straight lines. Their projections to the 
$q$-plane via $\pi$ are curves of degree $n$.
Some of them are shown in Figure \ref{fig:Tschirnhaus}.
\begin{remark}[Topology of the regularised energy surface]\quad\\ 
The  $n$--fold covering $\widehat{\Pi}$ 
of phase space $\widehat{P}$ has the form
\beq
(q,p) = \widehat{\Pi}(Q,P) = (Q^n,P\bar{Q}^{1-n}) \qquad 
\big((Q,P) \in \bC^\times \times \bC \big) \, ,
\Leq{loc:canonical}
with position $Q$ and momentum\,\footnote{We kept the mnemonic notation $P$ for the momentum on covering space, although $P$ also denotes the regularised phase space.} $P$. When restricted to $\bC^\times \times \bC$,
\[ \Phi: \bC^2 \to \bC^2  \qmbox{,} (Q,P)\mapsto \exp(\imath 2\pi/n) \, (Q,P) \]
generates the group of covering transformations $\Phi^k$ ($k\in \bZ/n\bZ$)
for $\widehat{\Pi}$.\\
The canonical symplectic form 
\[\omega := dq_1\wedge dp_1 + dq_2\wedge dp_2= \Re(dq\wedge d\bar{p})\] 
pulls back to 
$\Omega :=  \widehat{\Pi}^*(\omega) = n\,\Re(dQ\wedge d\bar{P})$.
So $\widehat{\Pi}$ is only locally canonical up to a factor $n$.
\begin{enumerate}[$\bullet$]
\item 
The lift of angular momentum 
\beq
L(q,p) = \mathrm{Im}(p\,\bar{q}) \qquad \big((q,p) \in \bC \times \bC \big)
\Leq{L:LC}
equals
\[ L \circ \widehat{\Pi} \, (Q,P) = \mathrm{Im}(P \, \bar{Q})  \qquad 
\big((Q,P) \in \bC^\times \times \bC \big) \, .\]
This uniquely extends to the real-analytic function 
\[\mathbf{L}: \bC\times \bC\to \bR\qmbox{,}\mathbf{L}(Q,P) 
=  \mathrm{Im}(P \, \bar{Q}) \, .\]
\item 
The lifted Hamiltonian equals 
\[ \textstyle
\widehat{H}\circ \widehat{\Pi} \, (Q,P) = \frac{\frac{|P|^2}{2m} -Z}{|Q|^{2(n-1)}}
\qquad
\big((Q,P) \in \bC^\times \times \bC \big)  \, .\]
Multiplying $\widehat{H}\circ \widehat{\Pi}  - E$ by $|Q|^{2(n-1)}$, 
this extends to the real-analytic function
\[ \textstyle
\frac{|P|^2}{2m} - E|Q|^{2(n-1)} - Z \qquad \big((Q,P) \in \bC\times \bC \big) \, .\] 
\begin{enumerate}[(a)]
\item 
For $E<0$ and $n\ge2$ its zero level is diffeomorphic to the standard sphere
$\bS^3\subseteq \bC\times \bC$, via radial projection.
Then the $\Phi^k$, $k=0,\ldots,n-1$,  
when restricted to the energy surface, are the covering transformations 
of \eqref{loc:canonical}, acting freely on the energy surface.
This implies that the regularised energy surface is diffeomorphic to the lens space
$L(n,1) = \bS^3/(\bZ/n\bZ)$.\\ 
In particular for the Kepler case we get real projective
space $\bR\mathrm{P}(3) = L(2,1)$. The fundamental group equals $\pi_1(L(n,1)) \cong
\bZ/n\bZ$.
\item 
For energy $E \ge 0$, for all $n \ge 2$ 
the regularised energy surface is diffeomorphic to  
$\bR^2\times \bS^1$, like for free motion in $d=2$ dimensions. \\
However, the non-trivial topology also leaves its traces here:
The Lagrange manifolds corresponding to
the scattering orbits with an arbitrary 
given initial direction fold over configuration space, with non-trivial
topological degree, see \cite[Section 5]{KK}.\hfill $\Diamond$
\end{enumerate}
\end{enumerate}
\end{remark}
%
\subsection{Proof of Theorem \ref{thm:complete}}\label{sub:proof}
We will begin with the definition of an open phase space neighbourhood 
$\widehat U^\vep \subseteq \widehat{P}$.
Then we will prove analyticity properties of phase space functions defined on it.
\subsubsection{The domain $\widehat U^\vep$}
The domain $\widehat U^\vep$ of a phase space chart, 
to be introduced in the following lemma, will be the intersection
of a domain $U^\vep \subseteq P$ of a chart on extended phase space $P$ 
and the domain  
$\widehat{P} \subseteq P$ of the $(q,p)$ chart, with 
$P = \widehat{P} \cup U^\vep$.
\begin{lemma}  \label{lem:U:vep}
We consider a maximal solution $(q,p): (t_-,t_+)\to\widehat U^\vep$
of the Hamiltonian equation for \eqref{Hamiltonian}, on
the open subset near the excluded fiber  $T^*_0\bR^d_q$
\beq
\widehat U^\vep := \left\{ (q,p) \in \widehat{P}\, \left|\ \|q\| < \vep\, ,\,
 \textstyle \widehat{H}(q,p) > \frac{-Z}{2n\|q\|^{\alpha_n}} \right.\right\} \qquad(\vep>0).
\Leq{U:epsilon}
Then $\frac{d}{dt}\LA q,p\RA  > \frac{Z}{n\vep^{\alpha_n}}$.
This implies that the length $ t_+ - t_->0$ of the time interval
has a uniform upper bound $2 \vep^{2-1/n} \sqrt{n/Z}$ and 
\begin{enumerate}[$\bullet$]
\item 
the solution either intersects the real-analytic hypersurface
\beq
\widehat{S}^\vep :=
\big\{ (q,p)\in \widehat{U}^\vep \; \big|\, \LA q,p\RA = 0\big\}
\Leq{hat:S:vep}
exactly once and transversally, in a {\bf pericenter} (meaning $\frac{d}{dt}\LA q,p\RA>0$), 
\item 
or it is a {\bf collision solution}, that is, $\lim_{t\searrow t_-} q(t) = 0\,$ or 
$\,\lim_{t\nearrow t_+} q(t)  =0$.
\item 
Every collision solution $(q,p):(t_-,t_+)\to\widehat{P}$ intersects $\widehat{U}^\vep$.
\end{enumerate}
\end{lemma}
\textbf{Proof:}
Note that $\widehat U^\vep$ is open in $T^*\bR^d$, too since 
$\widehat{P}\subseteq T^*\bR^d$ is open.
\begin{enumerate}[$\bullet$]
\item 
$\widehat{S}^\vep$ is a real-analytic hypersurface, since $q,p\neq0$ for
$(q,p)\in \widehat{U}^\vep$.
\item 
The squared distance from the origin is strictly convex as a function of time $t$:
\beq \textstyle
\eh \frac{d^2}{dt^2} \|q\|^2 = \frac{d}{dt} \LA q,p \RA 
= \frac{\|p\|^2}{m} -  \frac{\alpha_n Z}{\|q\|^{\alpha_n}}
= 2\widehat{H}(q,p) + \frac{2Z}{n\|q\|^{\alpha_n}} > 
 \frac{Z}{n\|q\|^{\alpha_n}} > \frac{Z}{n\vep^{\alpha_n}},
\Leq{convex}
using Definition \eqref{U:epsilon} in the last two inequalities.

So if $t\mapsto \|q(t)\|$ is strictly increasing on $(t_- , t_+)$, then 
$\lim_{t\searrow t_-} q(t)=0$ and
\beq
\vep^2 \ge \lim_{t\nearrow t_+}  \|q(t_+)\|^2 \ge 
\int_{t_-}^{t_+}(t-t_-) \textstyle \frac{2Z}{n \vep^{\alpha_n}}\, \mathrm{d}t
= (t_+ - t_-)^2\frac{Z}{n \vep^{\alpha_n}} \,.
\Leq{time:bound}
Equation \eqref{time:bound} leads to the upper bound
$\vep^{2-1/n} \sqrt{n/Z}$ for total time $t_+ - t_-$,
and $\lim_{t\searrow t_-} q(t)=0$ implies collision at $t_-$.\\
By reversibility of the flow, the strictly decreasing case is similar.
\item 
Otherwise there is a $t_0\in (t_-, t_+)$, with 
$\big( q(t_0) , p(t_0)\big)\in \widehat{S}^\vep$.
Then \eqref{time:bound} implies 
$\vep^2\ge  (t_\pm-t_0)^2\frac{Z}{ n \vep^{\alpha_n}}$, $q(t_0)$ is a pericenter, 
and the only time of intersection with $\widehat{S}^\vep$. In both cases we get 
the upper bound $2 \,\vep^{2-1/n} \sqrt{n/Z}$ for $t_+ - t_-$.
\item 
For a collision solution in $\widehat{P}$, by preservation of total energy $\widehat{H}$,
in the limit of collision time
$\|p\|^{2}/(2m) \sim \frac{Z}{\|q\|^{\alpha_n}}$.
But the condition $\widehat{H}(q,p) > \frac{-Z} {2n \|q\|^{\alpha_n}}$
in \eqref{U:epsilon} is equivalent to 
$\frac{\|p\|^{2}}{2m} > (1-\frac{1}{2n}) \frac{Z}{\|q\|^{\alpha_n}}$. \\
So the collision solution intersects~$\widehat{U}^\vep$.
\hfill $\Box$
\end{enumerate}
\begin{remarks}[Lemma \ref{lem:U:vep}] \label{rem:U:vep} \quad\\[-6mm]
\begin{enumerate}[1.]
\item 
{\bf (Definition of $\widehat U^\vep$)}
The factor $1/(2n)$ in the second condition of \eqref{U:epsilon} is not entirely arbitrary:
Only factors in the interval $(0,1/n)$ could be chosen, since
\begin{enumerate}[$\bullet$]
\item 
circular orbits have to be excluded from $\widehat U^\vep$, as they 
do not have a unique pericenter. But  for
them $\frac{\|p\|^{2}}{2m} = (1-\frac{1}{n})\frac{Z}{\|q\|^{\alpha_n}}$,  so that
$\widehat{H}(q,p) = \frac{-Z} {n \|q\|^{\alpha_n}}$.
\item 
Conversely all collision orbits should enter $\widehat U^\vep$, since we
want to regularise them by extending  $\widehat U^\vep$.
For them $\frac{\|p\|^{2}}{2m} \sim \frac{Z}{\|q\|^{\alpha_n}}$ in the limit of collision time.
\end{enumerate}
\item
{\bf (Repulsive case)}
For the repulsive case $Z<0$, 
where there are no collision orbits and no regularisation is needed, 
the second condition in \eqref{U:epsilon} is automatically fulfilled,
and the {\em estimates} of Lemma \ref{lem:U:vep} apply, too.
This may be useful for applications to $n$-particle systems interacting via 
Coulombic forces.
\item 
{\bf (Perturbations)}
Lemma \ref{lem:U:vep} persists under perturbations, with slightly worse constants:
Changing the potential in \eqref{Hamiltonian}
to $U(q):=\frac{Z}{\|q\|^{\alpha_n}}-W(q)$ with $W$ defined on a neighbourhood of the
origin and twice continuously differentiable,
\beq\textstyle
\frac{d}{dt}\!\LA q,p\RA
= \frac{\|p\|^2}{m} - \frac{\alpha_n Z}{\|q\|^{\alpha_n}} - \LA q,\nabla W(q)\RA 
> \frac{Z}{n \|q\|^{\alpha_n}} - \LA q,\nabla W(q)\RA
> \frac{Z}{2n\vep^{\alpha_n}},
\Leq{elapse}
for $\vep>0$ small, using \eqref{U:epsilon}. 
This would then substitute inequality \eqref{convex}.
\hfill $\Diamond$
\end{enumerate}
\end{remarks}
\subsubsection{Hamiltonian and angular momentum}\label{sub:sub:H:L}
We now regularise the motion by introducing new canonical coordinates
on $\widehat U^\vep$ and then extending $\widehat U^\vep$ to $U^\vep\subseteq P$.
These coordinates will be shown to be real-analytic and, taken together, will
rectify the flow. 

Total energy and angular momentum are treated first. They
are real-analytic from the outset.
We denote the restriction of the real-analytic {\bf Hamiltonian} \eqref{Hamiltonian} 
to $\widehat U^\vep$ by $\widehat H$, too.

{\bf Angular momentum} 
$\hat{L}:\widehat{U}^\vep\to \Lambda^2(\bR^d)$
equals $\hat{L}(q,p) = q\wedge p$. As a function with values in the Lie algebra
$\frak{so}(d,\bR) \cong \Lambda^2(\bR^d)$ it has entries 
\[\hat{L}_{i,j}(q,p) = q_j p_i - q_i p_j \, ,\]
and we set 
\beq
 \hat{\frak L}^{2}: \widehat U^\vep\to[0,\infty)\qmbox{,}  
 \hat{\frak L}^{2} := \eh\tr\big(\widehat{L}\widehat{L}^*\big) \, .
\Leq{frak:L}
Thus $\hat{\frak L}^{2}$ generalises the squared norm of the 
angular momentum vector from $d=3$ to arbitrary dimensions.
Both $\hat{L}$ and $\hat{\frak L}^{2}$ are real-analytic.

Depending on the value $E$ of the energy, we can bound the
image of $\hat{\frak L}^{2}$ in \eqref{frak:L} from above. 
If $\|q\| = \vep$, then $\hat{\frak L}^{2}(q,p)$ is maximal
if $\langle p,q\rangle =0$ and then takes the value 
$\eh \|p\|^2 \vep^2$.
In terms of $E=\frac{\|p\|^2}{2m}-U_n(\vep)$ we conclude that then
$\hat{\frak L}^{2}(q,p) = m(E+U_n(\vep))  \vep^2$. 
For $E<0$, $\|q\| = \vep$ may lie outside Hill's region,
but still the upper bound $\hat{\frak L}^{2}(q,p) \le m(E+U_n(\vep))  \vep^2$ is valid.
\subsubsection{Time as a phase space variable}\label{sub:sub:time}
Let $\widehat{T}: \widehat U^\vep \to\bR$ be the {\bf time} elapsed since the solution
was at its pericenter in $\widehat S^\vep$, respectively since collision. By (\ref{elapse}) 
there is only one such pericenter or collision of the orbit. 
So the phase space function $\widehat{T}$ is defined.
\begin{enumerate}[$\bullet$]
\item 
It is obvious that $\widehat{T}$ is {\bf real-analytic} 
on the open subset of $\widehat U^\vep$
where the angular momentum does not vanish. This follows from transversality
of the real-analytic flow to the real-analytic hypersurface $\hat{S}^\vep$.
\item 
To see that $\widehat{T}$ has this property everywhere, including on collision orbits,
we first consider dimension $d=2$. 
The lift of the domain $\widehat{U}^\vep$ defined in \eqref{U:epsilon}  
via $\widehat{\Pi}$ from \eqref{loc:canonical} equals
\[ \mathbf{\widehat{U}}^\vep := \widehat{\Pi}^{-1}(\widehat{U}^\vep) 
= \big\{ (Q,P)\in \bC^\times \times \bC\ \mid \|P\|^2 > 2m (1-\tfrac{1}{2n})Z \big\} \, .\] 
The defining equation is independent of $Q$, and we add the points with $Q=0$,
obtaining $\mathbf{U}^\vep \supseteq \mathbf{\widehat{U}}^\vep$. 
\[ \mathbf{\hat{S}}^\vep := \widehat{\Pi}^{-1}\big( \hat{S}^\vep \big)
\subseteq \mathbf{\widehat{U}}^\vep \] 
has the form $\Re(P\overline{Q}) = 0$, with $Q \neq 0$. 
As $P$ does not vanish on $\mathbf{U}^\vep$, we can drop the non-vanishing
condition for $Q$. We thus obtain a real-analytic hypersurface 
\[\mathbf{S}^\vep \subseteq  \mathbf{U}^\vep \, , \]
which is transversal with respect to the geodesic flow. 

Noticing that the factor $\frac{\mathrm{d}t}{\mathrm{d}s}$ in \eqref{time:rep} of 
time reparameterisation is real-analytic, by integrating $\frac{\mathrm{d}t}{\mathrm{d}s}$ 
along the orbit segment, we obtain a real-analytic function
$\mathbf{T} : \mathbf{U}^\vep \to \bR$. 
Its restriction to $\mathbf{\widehat{U}}^\vep$ is the lift 
\[\widehat{\mathbf{T}} = 
\widehat{T} \circ \widehat{\Pi} : \mathbf{\widehat{U}}^\vep \to\bR \]
of the above $\widehat{T}$. 
This shows that $\widehat{T}$ is real-analytic, since the Jacobi-Maupertuis
metric $\mathbf{g}_E$, see \eqref{lifted:metric}, is real-analytic in 
the position $Q$ and the parameter $E$.
\item
The generalisation to arbitrary dimensions $d\ge2$ is immediate, since  $\widehat{T}$
is invariant under the symplectic lift of the $\mathrm{SO}(d)$ action on $\widehat{M}$
to $\widehat{P}$.
\end{enumerate}
\begin{example}[Kepler case, $n=2$] \quad\\ 
The pericenter in $\widehat S^\vep$ has radius $\hat{r}_{\rm min}$, see 
Paragraph \ref{sub:sub:r:min} below.
There, 
Formula \eqref{rmin} includes the case of collision, with $\hat{L}=0$.\\
In the Kepler case, explicit integration is possible.
Then $\widehat{T}$ is given by
\beq
\widehat{T}(p,q) = \int^{\|q\|}_{{\hat{r}_{\rm min}}(p,q)}
\frac{\,\mathrm{d}r}{\sqrt{2\big(\widehat H(p,q) + Z r - \hat{\frak L}^2(p,q)/(2mr^{2})\big)}}\cdot
\sign(\LA q,p\RA) .
\Leq{peric:time}
Since $n$ is even,
the lower limit $\hat{r}_{\rm min}$  of the integral is real-analytic.
That $\widehat{T}$ is a real-analytic phase space function, too,
can be seen by explicit evaluation of the indefinite integral
\begin{eqnarray}
\lefteqn{\int \frac {r^{1/2}} {\sqrt{ 2Z + 2r E- \ell^2/(mr)}}\, \mathrm{d}r = }
\label{eq:explicit:int}\\ 
& &\hspace{-8mm}\textstyle
\frac{r}{\sqrt{2E}} \sqrt{1 + \frac{Z}{rE} - \frac{\ell^{2}}{2mr^{2}E}} -
              \frac{Z}{(2E)^{3/2}}
              \ln \left( Er + \eh Z + \sqrt{ E(r^{2}E + Zr - \ell^{2}/(2m)) }\right)
\nonumber
\end{eqnarray}
for $E>0$ (and similarly for  $E\le 0$, see {\sc Thirring} \cite[Chap.\ 4.2]{Th}):
With $\hat{r}_{\rm min}$ from \eqref{rmin} the contribution of the lower bound of integration
in \eqref{peric:time} simplifies to
\[\textstyle
\frac{-Z \log \left( \frac{ h l^2}{2m}+(Z/2)^2 \right)} {2 (2h)^{3/2}}.
\hfill \tag*{$\Diamond$} \]
\end{example}
\subsubsection{The generalised Laplace-Runge-Lenz vector}
Our last task is to construct a generalisation of the {\bf Laplace-Runge-Lenz vector}
of the Kepler case $n=2$. 
We will first work in $d=2$ dimensions and then generalise to arbitrary $d$.

\begin{figure}[h]
\centerline{
\includegraphics[height=40mm]{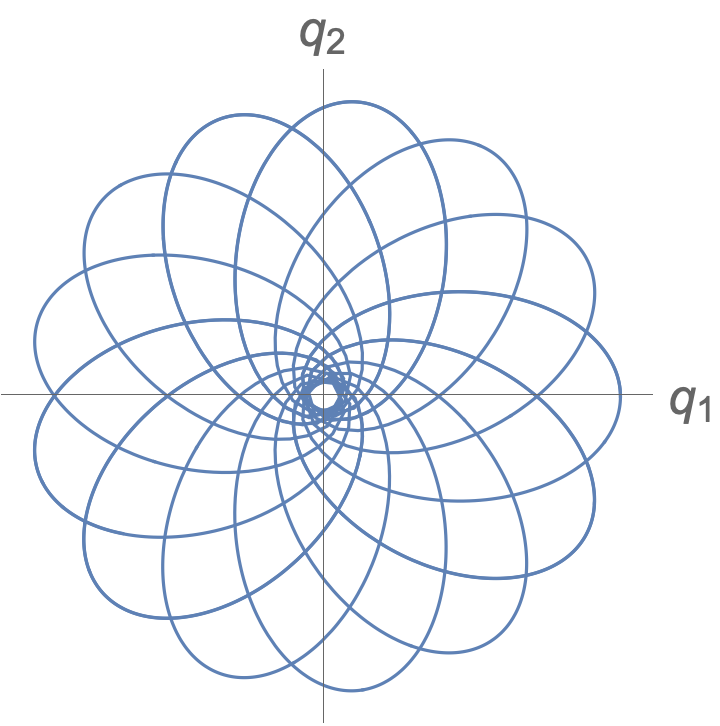}
\hspace{5mm}
\includegraphics[height=40mm]{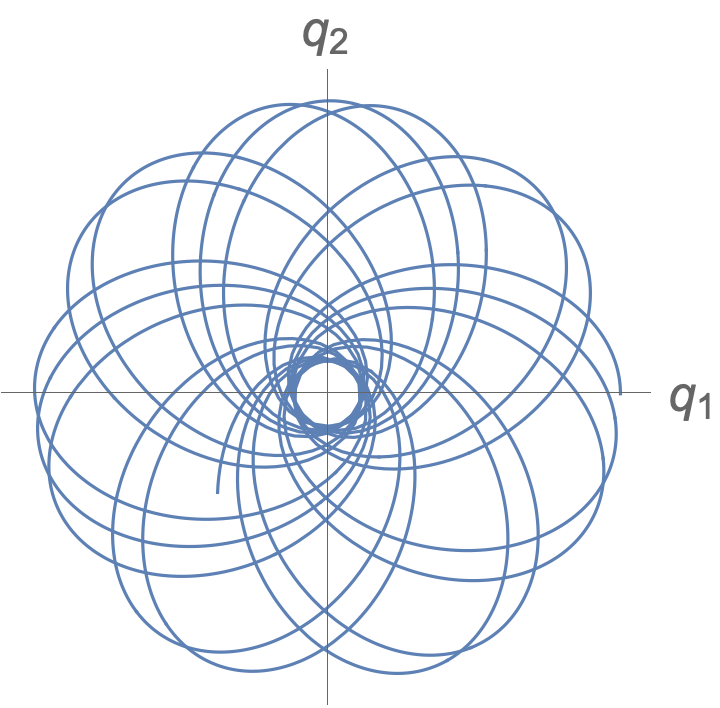}
\hspace{5mm}
\includegraphics[height=40mm]{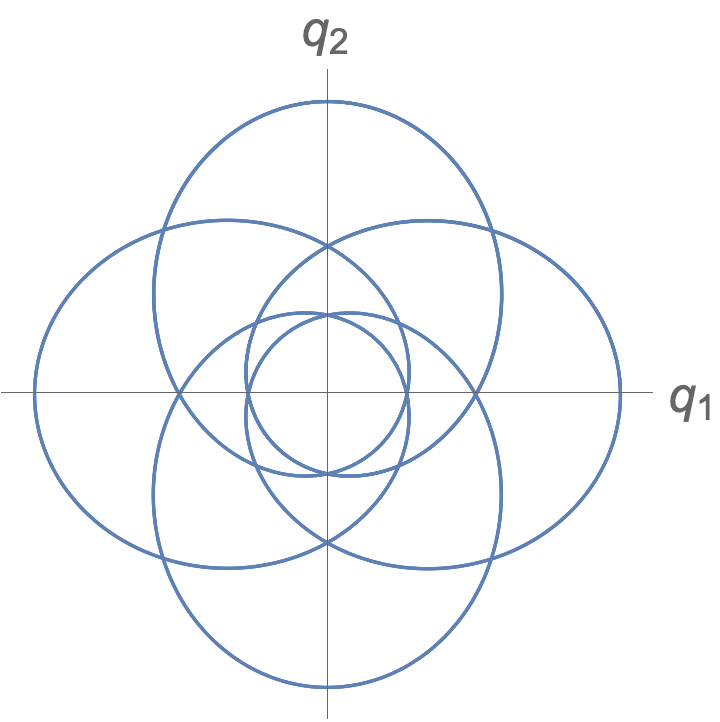}
}
\caption{Solutions for $n=3$ and negative energy. Left, right:
periodic orbit. Middle: quasiperiodic orbit. Angular momentum increases from left to right.}
\label{fig:neg:E}
\end{figure}

\begin{example}[Laplace-Runge-Lenz for the Kepler problem]\quad \label{ex:LRL}\\  
For $n=2$ this flow-invariant, real-analytic
function on $\widehat {P}$ has in complex coordinates the form
\beq
\textstyle
\widehat{V}:\widehat{P}\to \bC\qmbox{,}
\widehat{V}(q,p):= -\big( Z\frac{q}{|q|} + \imath \, p L(q,p) \big) \, .
\Leq{LRL}
Using \eqref{L:LC},
\beq
\textstyle \widehat{V}(q,p) 
= p \, \Re(p\bar{q}) - Z\frac{q}{|q|} = \eh p^2\bar{q} + q\widehat{H}(q,p) \, . 
\Leq{LRL:Ham}
So it vanishes exactly for the points on circular orbits.\\
Otherwise its direction $\widehat{V}/|\widehat{V}|$ is defined. For 
non-collision solutions it points towards the pericenter of the orbit, where 
$\imath \,p\bar{q}$ is real.

Note that for $n\ge 3$ a non-trivial 
generalisation of $\widehat{V}$ does not exist
as a function $\widehat{P}\to\bC$ on phase space, 
since then (following from Bertrand's Theorem for central-force potentials and 
analyticity) almost all orbits are not closed for negative energies, 
see Figure \ref{fig:neg:E}.
But here we are just to construct a map 
$\widehat{V}: \widehat{U}^\vep \to \bC$, and $ \widehat{U}^\vep$ does not contain 
any maximal orbit.
\hfill $\Diamond$
\end{example}
In order to locally generalise Example \ref{ex:LRL} to 
arbitrary $n\in\bN$, we start with the auxiliary function 
\[\frak{V}: \bC \times \bC \to \bC \qmbox{,}
\frak{V}(q,p) = p^n \bar{q}^{n-1} \, .\]
Via \eqref{loc:canonical}, this lifts to
\[ \frak{V} \circ \widehat{\Pi}(Q,P) = P^n\qquad
\big((Q,P)\in\bC^\times \times \bC \big),\]
which trivially extends to a real-analytic function on $\bC \times \bC$.
We now define the lifted Laplace-Runge-Lenz vector
\[\mathbf {V}:\mathbf{U}^\vep\to \bC \]
by setting $\mathbf {V}(Q,P)$ equal to $P_0^n$, where $(Q_0,P_0)\in \mathbf{S}^\vep$
is on the orbit through $(Q,P)$. 
By transversality of the flow with respect to $\mathbf{S}^\vep$, $\mathbf {V}$
is real-analytic.

For $n=2$, on $\mathbf{\widehat U}^\vep$ the argument 
$\mathbf {V} / |\mathbf{V}|$ of $\mathbf {V}$ coincides with
$\widehat{V} / |\widehat{V}| \circ \widehat{\Pi}$ in \eqref{LRL}.
As both are invariant under the respective flow, it suffices to compare them on $ \mathbf{S}^\vep$.
There, both terms in  \eqref{LRL:Ham} are dependent over the reals, and the first
one equals $\eh \frak{V}(q,p)$.

For $n\ge3$, too the generalised Laplace-Runge-Lenz vector 
transforms under orthogonal transformations via
\[ \widehat{V}(O q, O p) =O \widehat{V}(q,p)\qquad \big(O\in \mathrm{O}(2)\big),\]
since both $\frak{V}(q,p)$ and the flow transform covariantly under $\mathrm{O}(2)$.

For points $(q,p) \in \widehat{S}^\vep$ (which then belong to non-collision orbits), 
using \eqref{L:LC},
\[
\widehat{V}(q,p) =  \imath^{n-1} L(q,p)^{n-1} \, p \qmbox{and}
\widehat{V}(q,p) = \frac{ \imath^n L(q,p)^n}{\bar q}
= \frac{ \imath^n L(q,p)^n q}{|q|^2} \, . \] 
For $n$ even, this is a real multiple of $q$, for $n$ odd a real multiple of~$p$.

Collision orbits in $\bR^d$ do not lie in a unique two-plane, 
since $p$ and $q$ are linearly dependent. However,
$\widehat{V}(q,p)$ does not depend on the choice of the plane. 
Then the momentum 
 $P$ only depends on $p$ and $\|q\|$, and the direction of $P$ does does not
change with the flow. From \eqref{loc:canonical} we conclude that in these cases
$\widehat{V}(q,p)$ (which is still defined and non-zero via the covering construction) 
is linear dependent on both $p$ and $q$.

So in $d$ dimensions for all $n\in \bN$ the direction
\beq\textstyle
\widehat{A}: \widehat U^\vep \to \bS^{d-1}
\qmbox{,} \widehat{A} := \frac {\widehat{V}} {\|\widehat{V}\|}
\Leq{hat:A}
of the generalised Laplace-Runge-Lenz vector is well-defined and real-analytic.
Similarly,
\beq
\widehat{B}:\widehat U^\vep \to \bR^{d}
\qmbox{,}\widehat{B} :=\widehat{L} \widehat{A}
\Leq{hat:B}
is  real-analytic and perpendicular to $\widehat{A}$, by antisymmetry of $\widehat{L}$.
\subsubsection{Phase space functions that are {\em not} real-analytic} \label{sub:sub:r:min}
One might think that the pericentric radius is real-analytic, too, since this 
happens to be so for the
Kepler potential $U_2$. However, as we show now, this  only generalises to $U_n$ with
$n$ even.

$\hat{r}_{\rm min}:\widehat U^\vep \to[0,\vep)$ denotes 
the radius of the pericenter in $\widehat S^\vep$, and $\hat{r}_{\rm min}:=0$
for a collision orbit. It depends on the phase space point 
through $E$ and $\ell^2$:
$\hat{r}_{\rm min} = r_{\rm min}(\widehat H, \hat{\frak L}^{2})$.
In both cases $r_{\rm min}$ is given as a solution 
of the equation
\[ U_{\rm eff}(r) = -E\ \mbox{ , with }  
U_{\rm eff}(r) := U_n(r) -  \frac{\ell^{2} }{2r^2} \, . \]
Thus  $r_{\rm min}^2$  is a non-negative root of the polynomial 
\beq 
r^2 \mapsto(-Er^2+\eh \ell^2)^n -r^2 
\Leq{Polynomial}
of degree $n$, if $\alpha = \alpha_n$. There is exactly one such root if $E\ge0$
and two for $E<0$ and $\ell^2$ small. 
In the latter case $r_{\rm min}^2$ is the smaller one. 
We conclude that $r_{\rm min}^2$ is real-analytic in $E$ and $\ell^2$.
$r_{\rm min}$ is then real-analytic in $E$ and $\ell$ iff $n$ is even.
\begin{example} [Kepler case, $n=2$] 
Here the roots of \eqref{Polynomial} lead to
\beq
\hat{r}_{\rm min}(p,q) := 
\left\{ \begin{array}{cc}
     \frac{-Z+\sqrt{Z^2 + 2\widehat H(p,q)\hat{\frak L}^{2}(p,q)/m}}{2\widehat H(p,q)}
         &, \widehat H \neq 0 \\
     \eh\hat{\frak L}^{2}(p,q)/(m Z)  & , \widehat H = 0
     \end{array} \right. .
\Leq{rmin}
This is a real-analytic function on $\widehat U^\vep$. 
\hfill $\Diamond$
\end{example}
\subsubsection{Phase space extension}
We recapitulate what we have proven in Subsection \ref{sub:sub:H:L},
\ref{sub:sub:time}, in \eqref{hat:A} and in \eqref{hat:B}
for dimension $d\in \bN$ and $n\in\bN$:
There is a real-analytic map 
\beq
\widehat \Psi := \big(\widehat{T},\widehat{H}; \widehat{B},\widehat{A}\big):
\widehat U^\vep \ \longrightarrow \ T^*(\bR\times \bS^{d-1}) \, ,
\Leq{hat:psi}
and all entries except $\widehat{T}$ are constants of the motion.

The image $\widehat \Psi (\widehat U^\vep)$
misses the zero section $\cong \bR\times \bS^{d-1}$ of the
cotangent bundle, since the collision orbits are characterised by
zero angular momentum and thus $\widehat{B}=0$, 
and the point of collision on such an orbit corresponds to  $\widehat{T}\to 0$.

We now complete phase space by defining it as a set as the disjoint union
\beq
P:=\widehat{P} \ \sqcup \ (\bR\times \bS^{d-1})
\Leq{P:hatP:Kepler}
of a $2d$-dimensional and a $d$-dimensional manifold. We
think of $P$ as arising from $T^{*}\bR^d$ by exchanging the fiber 
$T^{*}_0\bR^d\cong \bR^d$ over $0\in \bR^d$ by $\bR\times \bS^{d-1}$.
\begin{remarks}\label{rem:one} \quad\\[-6mm]
\begin{enumerate}[1.]
\item {\bf (One spatial dimension)}
For $d=1$ we glue in $|\bS^0| = 2$ copies of the energy axis.
For given total energy $E$ we then have to choose how to glue two points 
with two ($E<0$) respectively four ($E\ge0$) open intervals to obtain a completion of
the energy surface. 
For one choice the particle changes between positions $q>0$ and $q<0$,
for the other it keeps on one side. 

In accordance with the unique completions for $d\ge2$, the first choice is natural 
for $n$ odd, the second for $n$ even.
In the latter case, including Kepler, in dimension $d=1$ 
the phase space $P$ is disconnected, unlike for $d\ge2$.
Then $P$ is homeomorphic to the disjoint union $\bR\times \bS^1\times \bS^0$
of two cylinders, for the particles keeping on the left respectively right half axis 
of configuration space.
\item 
{\bf (Dependence on $\vep$)}
The phase space region $\widehat U^\vep$ defined in \eqref{U:epsilon}
projects to a ball of radius $\vep$ in configuration space. 
Whereas some estimates of Lemma \ref{lem:U:vep} depend
on the value of $\vep$, the basic constructions leading to phase space regularisation 
do not. Consequently, several maps whose domain or image depend on
that value are not indexed by $\vep$. It will only be in the application of
Kepler regularisation, that the value of $\vep$ has to be chosen small
enough that the perturbative technique of Remark \ref{rem:U:vep}.3 is applicable. 
\hfill $\Diamond$
\end{enumerate}
\end{remarks}
$P$ becomes a real-analytic manifold by introducing a map $\Psi$ with domain
$U^\vep \subseteq P$,
\begin{align}
U^\vep &:=
\widehat U^\vep \sqcup (\bR\times \bS^{d-1}) 
\label{Uvep}\\
\Psi &:= \big(T,H;B,A\big): U^\vep \longrightarrow  T^*(\bR\times S^{d-1}) \, ,
\label{def:Psi}\\
\Psi\rstr_{ \widehat U^\vep} &:=  \widehat \Psi
\qmbox{and}
\Psi\rstr_{ \bR\times S^{d-1}}(h,a) := (0,h;0,a) \,. \nonumber
\end{align}
We define a topology on $U^\vep$ by declaring $\Psi$ to be a 
homeomorphism on its image.
and a real-analytic structure by the one inherited from $T^*(\bR\times \bS^{d-1})$.

Finally, we extend the pericentral hypersurface $\widehat{S}^\vep$ 
of $\widehat{U}^\vep$ defined in \eqref{hat:S:vep} to  
\beq
S^\vep := T^{-1}(0)\quad \mbox{, so that }\ 
S^\vep = \widehat{S}^\vep\sqcup (\bR\times \bS^{d-1}) \, . 
\Leq{S:vep}
As zero is a regular value of $T$, $S^\vep$ is a real-analytic hypersurface in $U^\vep$.
\subsubsection{The symplectic structure on $U^\vep$}
This subsection finishes the proof of Theorem \ref{thm:complete}, by showing that the
extension $\Psi$ of $\widehat \Psi$ is a diffeomorphism onto its image, and that the
corresponding extension of $\widehat \omega$ leads to a (known) symplectic structure.

\begin{lemma}[symplectic structure]  \label{lem:Poisson}
For all $n\in \bN$, $\widehat \Psi$ is a diffeomorphism onto its image. 
On the completed phase space region $U^\vep$ for all $1\le i,j\le d$
\[\{H,T\} = 1\qmbox{,}
\{H,A_i\} = \{H,B_i\} = 0\qmbox{,}
\{T,A_i\} = \{T,B_i\} = 0 \, , \]
\[
\{A_i,A_j\} =0 \qmbox{,}
\{A_i,B_j\} = \delta_{i,j}-A_iA_j\qmbox{,}
\{B_i,B_j\} = L_{i,j}  \, . \]
So $(P,\omega, H)$, with $\omega\rstr_{T^{*}(\bR^d\setminus \{0\})}$ 
the natural symplectic form and $H\in C^\omega (P,\bR)$, 
becomes a real-analytic complete Hamiltonian system, extending 
$(\widehat{P},\widehat{\omega}, \widehat{H})$. 
\end{lemma}
\begin{remark}[free and Kepler case] \label{rem:n:1:2}\quad\\
For $n=1$ and $n=2$ all these phase space functions are explicitly known and one
can calculate the Poisson brackets directly. For $n=2$, 
from the commutation relation 
$\{\widehat{V}_i,\widehat{V}_j\}=+2m\widehat{H}\widehat{L}_{i,j}$ 
of the Laplace-Runge-Lenz vector 
\eqref{LRL} (see, {\em e.g.}\ {\sc Thirring} \cite{Th}, \S4.2) we conclude that
the matrix of Poisson brackets equals
\[ \Big(\{\widehat{A}_i,\widehat{A}_j\}\Big)_{i,j}
= \frac{2m\widehat{H}} {\|\widehat{V}\|^4}
\Big(\|\widehat{V}\|^2 \widehat{L} - \widehat{V} \widehat{V}^\top \widehat{L}
- \widehat{L} \widehat{V} \widehat{V}^\top \Big)  \, . \]
So we show that 
$\widehat{V}\widehat{V}^\top \widehat{L} + \widehat{L}\widehat{V}\widehat{V}^\top
= \| \widehat{V} \|^2 \widehat{L}$ 
to prove that these Poisson brackets vanish.
For $d = 1$ both sides vanish.
The property ${\rm rank}\big(\widehat{L}(x)\big)\leq 2$ allows us to reduce
the remaining cases to dimension $d=2$, where the statement is obvious.\\[1mm]
As the flow generated by the Hamiltonian vector field $X_H$ is not explicitly given for
$n\ge3$, in these cases one needs abstract arguments to prove Lemma \ref{lem:Poisson}. 
\hfill $\Diamond$
\end{remark}
\textbf{Proof of Lemma \ref{lem:Poisson}:}\\[-6mm]
\begin{enumerate}[$\bullet$]
\item 
We noticed that $\widehat \Psi$ from \eqref{hat:psi} is real-analytic.
As a by-product of the calculation of the
Poisson brackets, it is a local diffeomorphism. $\widehat \Psi$ is injective, since 
the vectors $\widehat{A}$ and $\widehat{B}$  determine the plane or line of motion in
configuration space, and then energy, time to and direction of the pericenter
together with (scalar) angular momentum allow to reconstruct the particle position
on the trajectory. So $\widehat \Psi$ is a diffeomorphism onto its image.
\item 
Along solutions of the Hamiltonian equations $ \widehat T$ is their time parameter, 
up to an additive constant. So $\{\widehat{H},\widehat{T}\}=1$.\\
As $\widehat{A}$ and $\widehat{B}$ are constants of the motion, their 
Poisson brackets with $\widehat{H}$ vanish. 
\item 
To show that $\{\widehat{A}_i,\widehat{A}_j\} = 0$, like in Remark \ref{rem:n:1:2} 
we need only consider the case of $d=2$ dimensions. 
Then $\widehat{A\,}_1^2 + \widehat{A\,}_2^2 =1$, so that
$\{ \widehat{A}_1 , \widehat{A}_2 \} = 0$.
\item 
The Poisson brackets $\{\widehat{T},\widehat{A}_i\}$ vanish for the 
following reason: The Hamiltonian vector field of the function 
$(q,p)\mapsto \LA p,q\RA$
is tangential to the hypersurface 
$\widehat{S}^\vep=\widehat{T}^{-1}(0)$. It acts as a dilation and thus does
not change the vector $\widehat{A}$ on $\widehat{S}^\vep$. It follows that
$\{\widehat{T},\widehat{A}_i\}$ vanishes on $\widehat{S}^\vep$.
On the other hand 
$\frac{d}{dt}\{\widehat{T},\widehat{A}_i\}
=\big\{\widehat{H},\{\widehat{T},\widehat{A}_i\} \big\} = 0$, 
using the Jacobi identity, $\{\widehat{H},\widehat{T}\}=1$ and 
$\{\widehat{H},\widehat{A}_i\}=0$.
\item 
Since $\widehat{A}$ transforms like a vector, 
$\{\widehat{L}_{r,m},\widehat{A}_n\}
= \delta_{m,n}\widehat{A}_r - \delta_{r,n}\widehat{A}_m $.
So the Poisson brackets $\{\widehat{A}_i,\widehat{B}_j\}$ equal 
$\delta_{i,j} - \widehat{A}_i\widehat{A}_j$. 
\item 
The commutation relations for angular momentum are
\beq
 \big\{\widehat{L}_{i,j},\widehat{L}_{k,\ell} \big\} =
\delta_{i,\ell}\widehat{L}_{j,k}-\delta_{i,k}\widehat{L}_{j,\ell}
-\delta_{j,\ell}\widehat{L}_{i,k}+\delta_{j,k}\widehat{L}_{i,\ell} \, .
\Leq{LLL}
To show that
$\{\widehat{B}_i,\widehat{B}_j\}= \widehat{L}_{i,j}$, one substitutes the definition
\eqref{hat:B} of $\widehat{B}$, uses \eqref{LLL} and $\|\widehat{A}\|=1$. 

By real-analyticity of the variables and the symplectic form, 
and since $\widehat U^\vep$ is dense in $U^\vep$,
we can omit the hats in the Poisson brackets.
\item 
The flow of the Hamiltonian system
$(P,\omega, H)$ is complete, since the only incomplete trajectories of 
$(\widehat{P},\widehat{\omega}, \widehat{H})$ are the collision trajectories.
\hfill $\Box$
\end{enumerate}
\begin{remarks}\quad\label{rem:ps:tc:top}\\[-6mm] 
\begin{enumerate}[1.]
\item 
{\bf (Poisson structure)}
Using the Euclidean Riemannian metric on $T^*\bR^d\cong \bR^d_q \times \bR^d_p$, 
$T^* \bS^{d-1}$ embeds into $T^*\bR^d$. $T^*\bS^{d-1}$ contains 
the image of $(B,A)$ in \eqref{def:Psi} is a symplectic submanifold, 
given as $F^{-1}(0)$ for the regular value $0$ of
\[F\equiv(F_1,F_2):T^*\bR^d\to \bR^2 \qmbox{,}
(p,q) \mapsto \bigl( \LA q,q \RA - 1,\LA q,p \RA \bigr).\]
The corresponding Poisson structure given by Lemma \ref{lem:Poisson}
is the one of the Dirac bracket 
(see, {\em e.g.}\ {\sc Marsden} and {\sc Ratiu} \cite[Proposition 8.5.1]{MR}).
Namely, as $\{F_1,F_2\} =\LA q,q\RA$ equals $2$ on $T^*\bS^{d-1}$, 
restricted to $T^*\bS^{d-1}$ we have
\[ \{a,b\}_{\mathbf{Dirac}} = \{a,b\} + \eh \{a,F_1\}\{F_2,b\} - \eh\{a,F_2\}\{F_1,b\} \, . \]
for $a,b\in C^\infty(T^*\bR^d,\bR)$. 
So for indices $i,k = 1,\ldots,d$,
\[ \{q_i,q_k\}_{\mathbf{Dirac}} = 0 \mbox{ , }\ 
\{q_i,p_k\}_{\mathbf{Dirac}} = \delta_{i,k} - q_i q_k 
\mbox{ , }\ \{p_i,p_k\}_{\mathbf{Dirac}} = q_i p_k - q_k p_i \, . \]

In physics, these Poisson brackets are known as the ones of 
{\em Snyder} space-time \cite{Sn} (that is, de-Sitter space-time as a Lorentz manifold),
for the restrictions of momentum respectively position variables on phase space
$T^*\bS^{d-1}\subseteq T^*\bR^d$, see 
{\sc Leiva, Saavedra} and {\sc  Villanueva} \cite{LSV}, Section III.
\item
{\bf (Time change)} 
Although $(P, \omega,H)$ is a real-analytic Hamiltonian system, 
the projection from phase space $P$ to configuration space~$\bR^d_q$, 
\beq
\mathrm{Pr}: P\to \bR^d_q\mbox{ , }\ 
\mathrm{Pr} \big(T^*_q( \bR^d \setminus\{0\} )\big)  = \{q\} \mbox{ and } 
\mathrm{Pr} \big( \bR\times \bS^{d-1} \big) = \{0\}
\Leq{def:Pi:kepler}
is not smooth for $n>1$.
This follows, since the radius $\|q(t)\|$
of collision orbits as a function of time $t$ asymptotically varies like 
$t^{1/(2-1/n)}$,
and since time $T$ is one of the phase space coordinates on $U^\vep$.
\item 
{\bf (Topology of regularised phase space)}\\
Except for $d=1$ degrees of freedom the regularised phase space $(P,\omega)$
of the Kepler problem is not a cotangent bundle with its canonical symplectic form.

Consider $d=1$ first, see Remark \ref{rem:one}. 
There $P$ consists of two components homeomorphic to cylinders $\bS^1\times \bR$. 
When there is no collision ($\sign(q)\neq0$), then the value
of the Runge-Lenz map $A:U^\vep \to \bS^0$ equals $-\sign(q)$, 
and this sign is continuously extended to the collision. 
We thus need only consider one of these two components,
and show it to be diffeomorphic to $\bS^1\times \bR$.
We use the auxiliary map $r\in C^\infty(\bR^+, \bR^+)$, $r(q) = \sqrt{2m (Z/q +q)}$.
Hence $H(r(q),q) = q$.
Then we connect the end points of each interval 
$[-r(q),r(q)] \times \{q\} \subseteq \widehat{P} \subseteq P$
by the curve $[2T^-(x_0),0] \to P$, $t \mapsto \Phi(t,x_0)$, with $x_0 := (r(q),q)$ for the 
Hamiltonian flow $\Phi$.
These subsets of $P$ are mutually disjoint, 
homeomorphic to $\bS^1$, parameterized by  $q \in \bR^+$,
and their union equals $P$, see Figure \ref{fig:1d}. Smoothening the map yields a diffeomorphism.
\begin{figure}[htbp]
\begin{center}
\includegraphics[width=7cm]{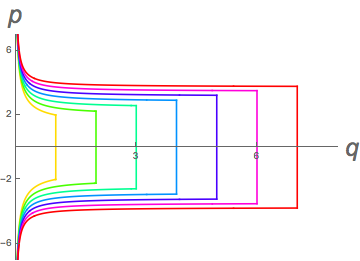}
\caption{Topology of regularised phase space in dimension $d=1$}
\label{fig:1d}
\end{center}
\end{figure}

For $d\in\bN$ and
energies $E<0$ the energy surface $H^{-1}(E)\subseteq P$ is 
compact.\footnote{By Moser regularisation it is homeomorphic to 
the unit tangent bundle $T^1 \bS^d$ of $\bS^d$.} This hypersurface
is the boundary of the two non-compact manifolds $H^{-1}((-\infty,E])$
and $H^{-1}([E,+\infty])$. For $d>1$ this contradicts the existence of 
a diffeomorphism $P\to T^*N$ for some $d$--dimensional manifold $N$. 
\hfill $\Diamond$
\end{enumerate}
\end{remarks}

\addcontentsline{toc}{section}{References}
\end{document}